\documentclass[12pt]{amsart}

\usepackage{amssymb,amsthm,mathtools,bm,color,scalerel}
\usepackage{hyperref}

\newtheorem{thm}{Theorem}
\newtheorem{oldthm}{Theorem}

\newtheorem{conj}{Conjecture}
\newtheorem{lemma}{Lemma}
\newtheorem{claim}[lemma]{Claim}
\newtheorem{prob}{Problem}
\newtheorem{cor}{Corollary}
\theoremstyle{definition}
\newtheorem{rem}{Remark}

\usepackage{tikz}
\usetikzlibrary{positioning}
\usetikzlibrary{svg.path}
\definecolor{orcid_color}{HTML}{A6CE39}

\hypersetup{pdfborder={0 0 0}} 
\DeclareRobustCommand{\orcidicon}{%
	\raisebox{.2mm}{\scalerel*{%
	\begin{tikzpicture}[xscale=1,yscale=-1,transform shape]
	\filldraw[color=orcid_color] svg {M256,128c0,70.7-57.3,128-128,128C57.3,256,0,198.7,0,128C0,57.3,57.3,0,128,0C198.7,0,256,57.3,256,128z};
	\filldraw[color=white] svg {M86.3,186.2H70.9V79.1h15.4v48.4V186.2z} svg {M108.9,79.1h41.6c39.6,0,57,28.3,57,53.6c0,27.5-21.5,53.6-56.8,53.6h-41.8V79.1z M124.3,172.4h24.5
		c34.9,0,42.9-26.5,42.9-39.7c0-21.5-13.7-39.7-43.7-39.7h-23.7V172.4z} svg {M88.7,56.8c0,5.5-4.5,10.1-10.1,10.1c-5.6,0-10.1-4.6-10.1-10.1c0-5.6,4.5-10.1,10.1-10.1
		C84.2,46.7,88.7,51.3,88.7,56.8z};
	\end{tikzpicture}}{|}}%
}
\newcommand{\orcid}[1]{\href{https://orcid.org/#1}{\orcidicon}}

\newcommand{\arxiv}[1]{arXiv:\href{https://doi.org/10.48550/arXiv.#1}{#1}}

\newcommand{\qbinom}[2]{\genfrac{[}{]}{0pt}{}{#1}{#2}}

\newcommand{\ksubsp}[2]{\Omega_{#1}^{\smash{(#2)}\vphantom{#2}}}

\title{Extremal problems for intersecting families of subspaces with 
a measure}
\author[Hajime Tanaka]{Hajime Tanaka\,\orcid{0000-0002-5958-0375}}
\address{\href{https://www.math.is.tohoku.ac.jp/index.html}{Research Center for Pure and Applied Mathematics}, Graduate School of Information Sciences, Tohoku University, Sendai 980-8579, Japan}
\email{htanaka@tohoku.ac.jp}
\urladdr{https://hajimetanaka.org/}
\author[Norihide Tokushige]{Norihide Tokushige\,\orcid{0000-0002-9487-7545}}
\address{College of Education, University of the Ryukyus, Nishihara  903-0213, Japan}
\email{hide@edu.u-ryukyu.ac.jp}
\urladdr{http://www.cc.u-ryukyu.ac.jp/~hide/}

\begin{document}
\begin{abstract}
We introduce a measure for subspaces of a vector space over a $q$-element
field, and propose some extremal problems for intersecting families.
These are $q$-analogue of Erd\H{o}s--Ko--Rado type problems, and we answer 
some of the basic questions.
\end{abstract}

\maketitle

\hypersetup{pdfborder={0 0 1}} 

\section{Introduction}
The purpose of the paper is to introduce a measure on the set of 
subspaces of a vector space over a finite field and propose some extremal 
problems for intersecting families of subspaces.
We then show some Erd\H{o}s--Ko--Rado type results for vector spaces with this measure.

To motivate our problems, we begin with the Erd\H{o}s--Ko--Rado theorem and its measure version.
Let $n$ and $k$ be positive integers with $n\geqslant k$.
Let $X_n:=\{1,2,\ldots, n\}$, and let $2^{X_n}$ and $\binom{X_n}{k}$ denote the power set of $X_n$ and the set of
$k$-element subsets of $X_n$, respectively.
A family $U\subset 2^{X_n}$ of subsets is called intersecting if 
$x\cap y\neq\emptyset$ for all $x,y\in U$.

\begin{oldthm}[\cite{EKR}]\label{thmA}
Let $\frac{k}{n}\leqslant \frac{1}{2}$.
If a family $U\subset\binom{X_n}{k}$ is intersecting, then
\begin{equation*}
	|U|/\binom{n}{k}\leqslant\frac{k}{n}.
\end{equation*}
Moreover, if $|U|/\binom{n}{k}=\frac{k}{n}$ and if $\frac{k}{n}< \frac{1}{2}$, then there exists $i\in X_n$ such that
\begin{equation*}
	U=\left\{x\in \binom{X_n}{k} : i\in x\right\}.
\end{equation*}
\end{oldthm}

This result has a measure counterpart.
Let $p$ be a real number with $0<p<1$.
We define a $p$-biased measure $\tilde{\mu}_p:2^{2^{X_n}}\to[0,1]$ by
\begin{equation*}
	\tilde{\mu}_p(U):=\sum_{x\in U}p^{|x|}(1-p)^{n-|x|}
\end{equation*}
for $U\subset 2^{X_n}$.
This is a probability measure and
\begin{equation*}
	\tilde{\mu}_p(2^{X_n})=\sum_{k=0}^n\binom{n}{k}p^k(1-p)^{n-k}=\left(p+(1-p)\right)^n=1
\end{equation*}
by the binomial theorem.
If $p=\frac{1}{2}$, then $\tilde{\mu}_{1/2}(U)=|U|/2^n$.
In general, we can think of $\tilde{\mu}_p(U)$ as a weighted sum of the subsets in $U$. If $p$ is fixed and $n$ is sufficiently large, then most of the weight of $\tilde{\mu}_p(2^{X_n})$ comes from the weights of subsets of size $k$ close to $np$.

\begin{oldthm}[\cite{AK}]\label{thmB}
Let $p\leqslant\frac{1}{2}$.
If a family $U\subset 2^{X_n}$ is intersecting, then
\begin{equation*}
	\tilde{\mu}_p(U)\leqslant p.
\end{equation*}
Moreover, if $\tilde{\mu}_p(U)=p$ and if $p<\frac{1}{2}$, then there exists $i\in X_n$ such that
\begin{equation*}
	U=\left\{x\in 2^{X_n} : i\in x\right\}.
\end{equation*}
\end{oldthm}

We find a correspondence between Theorems~\ref{thmA} and \ref{thmB} by comparing them.
The Erd\H{o}s--Ko--Rado theorem for $k$-uniform $t$-intersecting families by Wilson \cite{Wilson1984} and its $\tilde{\mu}_p$-measure counterpart by Friedgut \cite{Friedgut2008C} are one of the examples of such correspondences that occur not only between results but also between proofs.
In many cases, knowing the $k$-uniform version helps in understanding the corresponding $\tilde{\mu}_p$-measure version, and vice versa, even though there is no known direct implication between a statement of a $k$-uniform version and a statement of its $\tilde{\mu}_p$-measure version.

Now, we switch to work on subspaces from subsets.
Throughout the paper, we fix a prime power $q$.
Let $\mathbb{F}_q$ be the
$q$-element field, and let $\mathbb{F}_q^n$ denote the $n$-dimensional vector space over $\mathbb{F}_q$.
Let $\Omega_n$ and $\ksubsp{n}{k}$ denote the set of 
all subspaces of $\mathbb{F}_q^n$ and the set of $k$-dimensional subspaces of $\mathbb{F}_q^n$, respectively.
Define $[n]:=\frac{q^n-1}{q-1}$, $[n]!:=\prod_{j=1}^n[j]$, and $\qbinom{n}{k}:=\frac{[n]!}{[k]![n-k]!}$.
Then, $|\ksubsp{n}{k}|=\qbinom{n}{k}$.
A family $U\subset\Omega_n$ of subspaces is called
intersecting if $x\cap y\neq 0$ for all $x,y\in U$.

\begin{oldthm}[\cite{Hsieh1975DM}]\label{thmC}
Let $\frac{k}{n}\leqslant\frac{1}{2}$.
If a family $U\subset \ksubsp{n}{k}$ is intersecting, then
\begin{equation*}
	|U|/\qbinom{n}{k}\leqslant\frac{[k]}{[n]}.
\end{equation*}
Moreover, if $|U|/\qbinom{n}{k}=\frac{[k]}{[n]}$ and if $\frac{k}{n}< \frac{1}{2}$, then there exists $y\in \ksubsp{n}{1}$ such that
\begin{equation*}
	U=\left\{x\in \ksubsp{n}{k} : y\subset x\right\}.
\end{equation*}
\end{oldthm}

\begin{table}
\begin{tabular}{c||c|c}
 &uniform version & measure version\\
\hline\hline
subsets& Theorem~\ref{thmA} & Theorem~\ref{thmB}\\
\hline
subspaces& Theorem~\ref{thmC} & ?
\end{tabular}
\bigskip
\caption{}\label{table 1}
\end{table}

So far, we have mentioned three Erd\H{o}s--Ko--Rado type results.
It seems natural to expect a result that is a $q$-analogue of Theorem~\ref{thmB}, and
at the same time, a measure version corresponding to Theorem~\ref{thmC}
(see Table~\ref{table 1}).
To find such a result, we first need to introduce a measure on $\Omega_n$.
Let $\sigma$ be a positive real number, and let
\begin{equation*}
	\phi_{\sigma,n}(k):=\frac{\sigma^{k}q^{\binom{k}{2}}}{(-\sigma;q)_n},
\end{equation*}
where
\begin{equation}\label{q-binom}
	(-\sigma;q)_n:=\prod_{j=0}^{n-1}(1+\sigma q^j)
	=\sum_{k=0}^n\qbinom{n}{k} \sigma^k q^{\binom{k}{2}}.
\end{equation}
The second identity in \eqref{q-binom} is known as the $q$-binomial theorem.
(This is an identity as polynomials in $\sigma$ and hence is valid when, e.g., $\sigma<0$.)
Then, define a $\sigma$-biased measure 
$\mu_{\sigma}:2^{\Omega_n}\to[0,1]$ by
\begin{equation*}
	\mu_{\sigma}(U):=\sum_{x\in U}\phi_{\sigma,n}(\dim x)
\end{equation*}
for $U\subset\Omega_n$.
This is a probability measure and
\begin{equation*}
	\mu_{\sigma}(\Omega_n)=\sum_{k=0}^n\qbinom{n}{k}\phi_{\sigma,n}(k)=1.
\end{equation*}
We note that
\begin{equation*}
	\lim_{q\to 1}\qbinom{n}{k}\phi_{\sigma,n}(k)=
	\binom{n}{k} \!\! \left(\frac{\sigma}{1+\sigma}\right)^{\!\!k} \!
	\left(1-\frac{\sigma}{1+\sigma}\right)^{\!\! n-k}.
\end{equation*}
This suggests that the measure $\mu_{\sigma}$ is a $q$-analogue of the measure
$\tilde{\mu}_p$, and the maximum $\mu_{\sigma}$-biased measure of an intersecting
family of subspaces (in a possible result corresponding to 
Theorem~\ref{thmB}) should be $\frac{\sigma}{1+\sigma}$, which plays the role
of $p$ when $q\to 1$.
Indeed, we have the following.
\begin{thm}\label{thm1}
Let 
\begin{equation*}
	\sigma\leqslant q^{-2\lfloor \frac{n-1}{2}\rfloor-1} = \begin{cases} 
		q^{-n}&\text{if $n$ is odd,}\\
		q^{-n+1}&\text{if $n$ is even.}
	\end{cases}
\end{equation*}
If a family $U\subset\Omega_n$ is intersecting, then 
\begin{equation*}
 \mu_{\sigma}(U)\leqslant\frac{\sigma}{1+\sigma}.
\end{equation*}
Moreover, if $\mu_{\sigma}(U)=\frac{\sigma}{1+\sigma}$ and if $\sigma< q^{-2\lfloor \frac{n-1}{2}\rfloor-1}$, then there exists $y\in \ksubsp{n}{1}$ such that
\begin{equation*}
	U=\{x\in\Omega_n : y\subset x\}.
\end{equation*}
\end{thm}
Write $a=-2\lfloor \frac{n-1}{2}\rfloor-1$ for brevity.
Then, that $\sigma\leqslant q^a$ is equivalent to $\frac{\sigma}{1+\sigma}\leqslant\frac{q^a}{1+q^a}$, and we have $\frac{q^a}{1+q^a}\rightarrow \frac{1}{2}$ when $q\rightarrow 1$ (irrespective of the actual value of $a$).
Our proof of Theorem \ref{thm1} allows us to take the limit $q\rightarrow 1$, and we thus restore Theorem \ref{thmB}.
However, unlike the bound $p\leqslant\frac{1}{2}$ in Theorem~\ref{thmB}, the bound $\sigma\leqslant q^a$ in Theorem \ref{thm1} is not best possible in general.
Recall that $\frac kn$ in Theorem \ref{thmA} corresponds to $p$ in Theorem \ref{thmB}.
Similarly, it seems natural to expect that $\frac{[k]}{[n]}$ in Theorem \ref{thmC} corresponds to $\frac{\sigma}{1+\sigma}$ in its measure version.
Thus, let $\frac{\sigma}{1+\sigma}=\frac{[pn]}{[n]}$, or equivalently, $\sigma=\frac{[pn]}{[n]-[pn]}$ for a fixed $p$ with $0<p<1$, where we extend the notation $[\lambda]=\frac{q^{\lambda}-1}{q-1}$ to any $\lambda\in\mathbb{R}$.
Then, we have $\frac\sigma{1+\sigma}=\frac{q^{pn}-1}{q^n-1}\to p$ as $q\to 1$.

\begin{conj}\label{conj1}
For every $0<p<\frac{1}{2}$, there exists $n_0$ such that 
the following holds for all $n>n_0$ and $\sigma=\frac{[pn]}{[n]-[pn]}$:
if a family $U\subset\Omega_n$ is intersecting, then
\begin{equation*}
	\mu_{\sigma}(U)\leqslant\frac{\sigma}{1+\sigma},
\end{equation*}
with equality if and only if there exists $y\in \ksubsp{n}{1}$ such that
\begin{equation*}
	U=\{x\in\Omega_n : y\subset x\}.
\end{equation*}
\end{conj}
Note that Conjecture~\ref{conj1} does not cover Theorem~\ref{thm1}.
Indeed, if $p$ is fixed and $n$ is sufficiently large, then 
\begin{equation}\label{eq1}
	\sigma=\frac{[pn]}{[n]-[pn]}\sim q^{-(1-p)n}.
\end{equation}
Thus, Conjecture~\ref{conj1} applies to the case roughly 
$q^{-n}<\sigma<q^{-\frac{n}{2}}$.
We were unable to prove (or disprove) this conjecture.
Instead, we present some weaker results supporting it under more general settings.

Let $t$ be a fixed positive integer.
A family $U\subset \Omega_n$ of subspaces is called $t$-intersecting if 
$\dim(x\cap y)\geqslant t$ for all $x,y\in U$.
Fix $y\in\ksubsp{n}{t}$, and define a
$t$-intersecting family $A_n^{(t)}$ by
\begin{equation}\label{def:An}
	A_n^{(t)}:=\{x\in\Omega_n:y\subset x\}.
\end{equation}
Note that $A_n^{(1)}$ is an optimal ($1$-)intersecting family in Theorem \ref{thm1} and Conjecture \ref{conj1}.
Using \eqref{q-binom} (with $\sigma$ replaced by $\sigma q^t$) and $\sigma^kq^{\binom{k}{2}}=\sigma^tq^{\binom{t}{2}}\cdot(\sigma q^t)^{k-t}q^{\binom{k-t}{2}}$, we have
\begin{equation*}
	\sum_{k=t}^n\qbinom{n-t}{k-t}\sigma^kq^{\binom{k}{2}}= \sigma^tq^{\binom{t}{2}}(-\sigma q^t;q)_{n-t},
\end{equation*}
and so
\begin{equation}\label{mu_s(A)}
	\mu_{\sigma}(A_n^{(t)})=\sum_{k=t}^n\qbinom{n-t}{k-t}\phi_{\sigma,n}(k)
	=\frac{\sigma^tq^{\binom{t}{2}}}{(-\sigma;q)_t}=\left(\!-\frac{1}{\sigma};\frac{1}{q}\right)_{\!t}^{\!\!-1}.
\end{equation}
In particular, $\mu_{\sigma}(A_n^{(1)})=\frac{\sigma}{1+\sigma}$.
We are interested in the maximum $\sigma$-biased measure of $t$-intersecting
families, and let
\begin{equation*}
	f(n,t,\sigma):=\max\{\mu_{\sigma}(U)^{\frac{1}{n}}:\text{$U\subset\Omega_n$ is 
	$t$-intersecting}\}.
\end{equation*}
\begin{prob}
Find a condition for $\sigma$ that guarantees
$f(n,t,\sigma)=\mu_{\sigma}(A_n^{(t)})^{\frac{1}{n}}$.
\end{prob}
Based on \eqref{eq1}, we define
\begin{equation*}
	\sigma_{\theta,n}:=q^{-(1-\theta)n},
\end{equation*}
and write
\begin{equation*}
	\mu_{\theta,n}:=\mu_{\sigma_{\theta,n}}.
\end{equation*}
Then, for $0<\theta<\frac{1}{2}$, we have
\begin{equation*}
	\lim_{n\to\infty}\mu_{\theta,n}(A_n^{(t)})^{\frac{1}{n}}=q^{-(1-\theta)t}
\end{equation*}
(see Lemma~\ref{la2} in Section \ref{sec3}), and this is the best we can do approximately as shown below.
\begin{thm}\label{thm2}
We have
\begin{equation*}
	\lim_{n\to\infty}  f(n,t,\sigma_{\theta,n})=
	\begin{cases}
		q^{-(1-\theta)t} & \text{if $0<\theta<\frac{1}{2}$,}\\
		1 & \text{if $\frac{1}{2}<\theta<1$.}
	\end{cases}
\end{equation*}
\end{thm}

\begin{conj}
We have
\begin{equation*}
	\lim_{n\to\infty}  f(n,t,\sigma_{\frac{1}{2},n})= q^{-\frac{1}{2}t}.
\end{equation*}
\end{conj}

Two families $U,W\subset\Omega_n$ are called cross $t$-intersecting
if $\dim(x\cap y)\geqslant t$ for all $x\in U$ and $y\in W$.
Let
\begin{equation*}
	g(n,t,\sigma_1,\sigma_2):=
	\max\{\left(\mu_{\sigma_1}(U)\mu_{\sigma_2}(W)\right)^{\frac{1}{n}}:
	\text{$U,W\subset\Omega_n$ are cross $t$-intersecting}\}.
\end{equation*}

\begin{thm}\label{thm3}
If $0<\theta_1,\theta_2<\frac{1}{2}$, then
$\displaystyle\lim_{n\to\infty}g(n,t,\sigma_{\theta_1,n},\sigma_{\theta_2,n})=
q^{-(2-\theta_1-\theta_2)t}$.
\end{thm}

If $U\subset\Omega_n$ is $t$-intersecting, then $U,U$ are cross
$t$-intersecting.
This gives us that
\begin{equation*}
	f(n,t,\sigma)^2 \leqslant g(n,t,\sigma,\sigma).
\end{equation*}
Thus, Theorem~\ref{thm3} implies Theorem~\ref{thm2} 
for the case $0<\theta<\frac{1}{2}$.

In Section~\ref{sec2}, we prove Theorem~\ref{thm1}.
To this end, we first translate the problem into a semidefinite programming problem and then
solve it by computing eigenvalues of related matrices.
In Section~\ref{sec3}, we prove Theorem~\ref{thm2} by a probabilistic approach.
For this, we use that the distribution $\qbinom{n}{k}\phi_{\sigma,n}(k)$ on the points $k=0,1,\dots,n$ is concentrated around $\theta n$.
We also use the result (Theorem~\ref{thmD}) about the 
maximum size of $t$-intersecting families of subspaces of dimension $k$ due to Frankl and Wilson \cite{FW1986}.
In a similar way, we prove Theorem~\ref{thm3} in Section~\ref{sec4}, where we need the result (Theorem \ref{thmG}) about cross $t$-intersecting uniform families due to Cao, Lu, Lv, and Wang \cite{CLLW}.
We mention that Theorem \ref{thmG} partly generalizes the result (Theorem~\ref{thmF}) about the case $t=1$ due to Suda and Tanaka \cite{ST2014BLMS}, which was also proved by solving the corresponding semidefinite programming problem.

\section{Proof of Theorem~\ref{thm1}}\label{sec2}

For a non-empty finite set $\Lambda$, let $\mathbb{R}^{\Lambda}$ be the set of real column vectors with coordinates indexed by $\Lambda$.
For two non-empty finite sets $\Lambda$ and $\Xi$, we also identify $\mathbb{R}^{\Lambda\times\Xi}$ with the set of real matrices with rows indexed by $\Lambda$ and columns indexed by $\Xi$.
When $\Lambda\subset\Lambda'$ and $\Xi\subset\Xi'$, we often view $\mathbb{R}^{\Lambda}$ (resp.~$\mathbb{R}^{\Lambda\times\Xi}$) as a subspace of $\mathbb{R}^{\Lambda'}$ (resp.~$\mathbb{R}^{\Lambda'\times\Xi'}$) in the obvious manner.
Define $W_{k,\ell},\overline{W}_{\!k,\ell}\in\mathbb{R}^{\ksubsp{n}{k}\times \ksubsp{n}{\ell}}$ by
\begin{equation*}
	(W_{k,\ell})_{x,y}=\begin{cases} 1 & \text{if}\ x\subset y \ \text{or} \ x\supset y, \\ 0 & \text{otherwise}, \end{cases} \qquad (\overline{W}_{\!k,\ell})_{x,y}= \begin{cases} 1 & \text{if}\ x\cap y=0, \\ 0 & \text{otherwise}, \end{cases}
\end{equation*}
for $x\in\ksubsp{n}{k}$, $y\in\ksubsp{n}{\ell}$.

We define the subspaces $\bm{U}_i$ $(0\leqslant i\leqslant \lfloor \frac{n}{2}\rfloor)$ of $\mathbb{R}^{\Omega_n}$ by
\begin{equation*}
	\bm{U}_i=\{\bm{u}\in\mathbb{R}^{\ksubsp{n}{i}}:W_{i-1,i}\bm{u}=0\} \qquad (0\leqslant i\leqslant \lfloor \tfrac{n}{2}\rfloor),
\end{equation*}
where $W_{-1,0}:=0$.
Since
\begin{equation*}
	W_{i-1,i}W_{i,i-1}=W_{i-1,i-2}W_{i-2,i-1}+q^{i-1}\qbinom{n-2i+2}{1}W_{i-1,i-1} \qquad (1\leqslant i\leqslant \lfloor \tfrac{n}{2}\rfloor)
\end{equation*}
and the RHS is positive definite, it follows that the matrices $W_{i-1,i}$ $(1\leqslant i\leqslant \lfloor \tfrac{n}{2}\rfloor)$ have full rank $\qbinom{n}{i-1}$ and hence
\begin{equation*}
	\dim \bm{U}_i=d_i:=\qbinom{n}{i}-\qbinom{n}{i-1} \qquad (0\leqslant i\leqslant \lfloor \tfrac{n}{2}\rfloor).
\end{equation*}

For $0\leqslant i\leqslant \lfloor \tfrac{n}{2}\rfloor$, we fix an orthonormal basis $\bm{u}_{i,1},\bm{u}_{i,2},\dots,\bm{u}_{i,d_i}$ of $\bm{U}_i$, and define
\begin{equation}\label{ONB}
	\bm{u}_{i,r}^k = q^{-\frac{i(k-i)}{2}} \qbinom{n-2i}{k-i}^{-\frac{1}{2}} W_{k,i} \bm{u}_{i,r} \qquad (1\leqslant r\leqslant d_i, \ i\leqslant k\leqslant n-i).
\end{equation}
In \cite{ST2014BLMS}, it is shown that the $\bm{u}_{i,r}^k$ form an orthonormal basis of $\mathbb{R}^{\Omega_n}$, and that
\begin{equation}\label{block-diagonalization}
	\overline{W}_{\!k,\ell}\bm{u}_{i,r}^{\ell} = \theta_i^{k,\ell} \bm{u}_{i,r}^k,
\end{equation}
where $\bm{u}_{i,r}^k:=0$ if $k<i$ or $k>n-i$, and
\begin{align*}
	\theta_i^{k,\ell} &= (-1)^i q^{\binom{i}{2}+k\ell-\frac{i(k+\ell)}{2}} \qbinom{n-k-i}{\ell-i} \qbinom{n-2i}{k-i}^{\frac{1}{2}} \qbinom{n-2i}{\ell-i}^{-\frac{1}{2}} \\
	&= (-1)^i \frac{ q^{\binom{i}{2}+k\ell-\frac{i(k+\ell)}{2}} }{ (q;q)_{n-k-\ell} } \!\left( \frac{ (q;q)_{n-k-i} (q;q)_{n-\ell-i} }{ (q;q)_{k-i} (q;q)_{\ell-i} } \right)^{\!\frac{1}{2}}.
\end{align*}
We note that
\begin{equation*}
	\theta_i^{k,\ell}=\theta_i^{\ell,k} \qquad (i\leqslant k,\ell\leqslant n-i).
\end{equation*}
See also \cite{FW1986}.
Thus, for $0\leqslant i\leqslant \lfloor \frac{n}{2}\rfloor$ and $1\leqslant r\leqslant d_i$, the subspace
\begin{equation*}
	\bm{V}_{i,r}=\operatorname{span}\{\bm{u}_{i,r}^i,\bm{u}_{i,r}^{i+1},\dots,\bm{u}_{i,r}^{n-i}\}
\end{equation*}
of $\mathbb{R}^{\Omega_n}$ is invariant under all the $\overline{W}_{\!k,\ell}$, and we have
\begin{equation*}
	\mathbb{R}^{\Omega_n}=\bigoplus_{i=0}^{\lfloor \frac{n}{2}\rfloor}\bigoplus_{r=1}^{d_i} \bm{V}_{i,r} \qquad \text{(orthogonal direct sum)}.
\end{equation*}
We have $d_0=1$, and without loss of generality, we set $\bm{u}_{0,1}$ to be the vector with $1$ in coordinate $0\in\ksubsp{n}{0}$, and $0$ in all other coordinates.

Let $\Delta\in\mathbb{R}^{\Omega_n\times\Omega_n}$ be the diagonal matrix with diagonal entries
\begin{equation*}
	\Delta_{x,x} =\mu_{\sigma}(x) = \phi_{\sigma,n}(\dim x) \qquad (x\in\Omega_n),
\end{equation*}
and let $J\in\mathbb{R}^{\Omega_n\times\Omega_n}$ be the matrix all of whose entries are $1$.
Let $S\mathbb{R}^{\Omega_n\times\Omega_n}$ denote the set of symmetric matrices in $\mathbb{R}^{\Omega_n\times\Omega_n}$.
Following \cite{Lovasz1979IEEE,Schrijver1979IEEE}, we formulate the problem of maximizing the $\mu_{\sigma}$-biased measure of an intersecting family into a semidefinite programming problem as follows:
\begin{equation*}
\begin{array}{lll}
	\text{(P):} & \text{maximize} & \operatorname{tr}(\Delta J\Delta X) \\[0.12in]
		& \text{subject to} & \operatorname{tr}(\Delta X) = 1, \ X \succcurlyeq 0, \ X \geqslant 0, \\
		&& X_{x,y} = 0 \ \ \text{for} \ x,y\in \Omega_n, \, x\cap y=0,
\end{array}
\end{equation*}
where $X \in S\mathbb{R}^{\Omega_n\times\Omega_n}$ is the variable, $\operatorname{tr}$ means trace, and $X \succcurlyeq 0$ (resp.~$X \geqslant 0$) means that $X$ is positive semidefinite (resp.~nonnegative).
Indeed, if $\bm{x}\in\mathbb{R}^{\Omega_n}$ is the characteristic vector of a non-empty intersecting family $U\subset\Omega_n$, then the matrix
\begin{equation}\label{X from U}
	X:=\mu_{\sigma}(U)^{-1}\bm{x}\bm{x}^{\mathsf{T}}\in S\mathbb{R}^{\Omega_n\times\Omega_n}
\end{equation}
satisfies all the constraints and we have $\operatorname{tr}(\Delta J\Delta X)=\mu_{\sigma}(U)$.
We recommend the introductory paper \cite{Todd2001AN} on semidefinite programming.
We note that the above semidefinite programming problem can be generalized to handle cross-intersecting families.
See \cite{ST2014BLMS,STT2017MP}, and also \cite{LPV2023pre}.
The dual problem for (P) is then given by
\begin{equation*}
\begin{array}{lll}
	\text{(D):} & \text{minimize} & \alpha \\
		& \text{subject to} & S:= \alpha \Delta -\Delta J \Delta +A -Z \succcurlyeq 0, \ Z \geqslant 0, \\
		&& A_{x,y} =0 \ \ \text{for} \ x,y\in \Omega_n, \, x\cap y\ne 0,
\end{array}
\end{equation*}
where $\alpha\in\mathbb{R}$ and $A,Z\in S\mathbb{R}^{\Omega_n\times\Omega_n}$ are the variables.
For any feasible solutions to (P) and (D), we have
\begin{equation}\label{weak duality}
	\alpha-\operatorname{tr}(\Delta J\Delta X) = \operatorname{tr}((\alpha\Delta-\Delta J\Delta)X) = \operatorname{tr}((S-A+Z)X) \geqslant 0,
\end{equation}
since $\operatorname{tr}(SX)\geqslant 0$, $\operatorname{tr}(ZX)\geqslant 0$, and $\operatorname{tr}(AX)=0$, and hence $\alpha$ gives an upper bound on $\mu_{\sigma}(U)$.

Our goal now is to find a feasible solution $(\alpha,A,Z)$ to (D) with
\begin{equation*}
	\alpha:=\frac{\sigma}{1+\sigma}.
\end{equation*}
Instead of working directly with the matrix $S$ above, we consider the positive semidefiniteness of
\begin{equation*}
	S':=\Delta^{-\frac{1}{2}}S\Delta^{-\frac{1}{2}}=\alpha I- \Delta^{\frac{1}{2}}J\Delta^{\frac{1}{2}} +A' -Z',
\end{equation*}
where $Z'\geqslant 0$ and $(A')_{x,y}=0$ whenever $x\cap y\ne 0$.
We set
\begin{equation*}
	Z':=0,
\end{equation*}
and choose $A'$ of the form
\begin{equation}\label{A'}
	A'=\sum_{k+\ell\leqslant n} \frac{a_{k,\ell}'}{\theta_0^{k,\ell}} \overline{W}_{\!k,\ell},
\end{equation}
where $a'_{k,\ell}\in\mathbb{R}$.
Let $\bm{1}\in\mathbb{R}^{\Omega_n}$ be the all-ones vector and recall our choice of the vector $\bm{u}_{0,1}\in\bm{U}_0=\mathbb{R}^{\ksubsp{n}{0}}$.
Then, we have (cf.~\eqref{ONB})
\begin{equation*}
	\bm{1}=\sum_{k=0}^n W_{k,0}\bm{u}_{0,1} =\sum_{k=0}^n \qbinom{n}{k}^{\frac{1}{2}} \bm{u}_{0,1}^k,
\end{equation*}
so that
\begin{equation}\label{special eigenvector}
	\Delta^{\frac{1}{2}}\bm{1}= \sum_{k=0}^n \phi_k^{\frac{1}{2}} \qbinom{n}{k}^{\frac{1}{2}} \bm{u}_{0,1}^k,
\end{equation}
where we abbreviate
\begin{equation*}
	\phi_k:=\phi_{\sigma,n}(k) \qquad (0\leqslant k\leqslant n).
\end{equation*}

Let $\mathbb{R}^{(n+1)\times(n+1)}$ be the set of real matrices with rows and columns indexed by $0,1,\dots,n$.
Let $C\in\mathbb{R}^{(n+1)\times(n+1)}$ be the lower triangular matrix defined by
\begin{equation*}
	C_{k,\ell}=\qbinom{k}{\ell} \qquad (0\leqslant \ell\leqslant k\leqslant n).
\end{equation*}
Then, it follows that
\begin{equation*}
	(C^{-1})_{k,\ell}=\qbinom{k}{\ell}(-1)^{k-\ell}q^{\binom{k-\ell}{2}} \qquad (0\leqslant \ell\leqslant k\leqslant n).
\end{equation*}
Indeed, for $0\leqslant \ell\leqslant k\leqslant n$, we have
\begin{equation*}
	\sum_{j=\ell}^k C_{k,j} \qbinom{j}{\ell}(-1)^{j-\ell}q^{\binom{j-\ell}{2}} = \qbinom{k}{\ell}\sum_{j=\ell}^k \qbinom{k-\ell}{j-\ell}(-1)^{j-\ell}q^{\binom{j-\ell}{2}} = \qbinom{k}{\ell}(1;q)_{k-\ell} =\delta_{k,\ell}
\end{equation*}
by \eqref{q-binom} (with $\sigma=-1$).
Now, let $G\in\mathbb{R}^{(n+1)\times(n+1)}$ be the upper triangular matrix given by
\begin{equation*}
	G_{k,\ell}=\qbinom{n-k}{n-\ell} (-1)^k\sigma^{\ell}q^{\binom{k}{2}+\binom{\ell}{2}} \qquad (0\leqslant k\leqslant\ell\leqslant n),
\end{equation*}
and let
\begin{equation}\label{G and F}
	F=CGC^{-1}.
\end{equation}
Note that the diagonal entries of $G$ are
\begin{equation*}
	(-1)^k\sigma^kq^{k(k-1)} \qquad (0\leqslant k\leqslant n),
\end{equation*}
and these are the eigenvalues of $G$, and hence of $F$.
We will set
\begin{equation}\label{our choice for A'}
	a_{k,\ell}' := \frac{F_{k,\ell}}{1+\sigma}\cdot\frac{\phi_k^{\frac{1}{2}}\qbinom{n}{k}^{\frac{1}{2}}}{\phi_{\ell}^{\frac{1}{2}}\qbinom{n}{\ell}^{\frac{1}{2}}} \qquad (0\leqslant k\leqslant n,\, 0\leqslant\ell\leqslant n-k)
\end{equation}
in \eqref{A'}, and show that the corresponding $S$ gives an optimal feasible solution to (D).

To describe the matrix $F$, we will use the following lemma.

\begin{lemma}\label{technical identity}
For integers $a,b$, and $c$ such that $a\geqslant b\geqslant 0$ and $c\geqslant 0$, we have
\begin{equation*}
	\sum_{j=0}^b \qbinom{b}{j}\qbinom{a-j}{c}(-1)^jq^{\binom{j}{2}}=q^{b(a-c)} \qbinom{a-b}{c-b}.
\end{equation*}
In particular, the LHS above vanishes if $b>c$.
\end{lemma}

\begin{proof}
First, we have
\begin{equation*}
	\qbinom{a-j}{c} = \sum_{d=0}^c q^{d(a-b-c+d)} \qbinom{a-b}{c-d} \qbinom{b-j}{d} \qquad (0\leqslant j\leqslant b).
\end{equation*}
To see this, fix $z\in\ksubsp{a-j}{a-b}$ and count $x\in\ksubsp{a-j}{c}$ such that $\dim (x\cap z)=c-d$ for each $d$ $(0\leqslant d\leqslant c)$.
Then, using \eqref{q-binom} (with $\sigma=-1$), we have
\begin{align*}
	\sum_{j=0}^b \qbinom{b}{j}\qbinom{a-j}{c}(-1)^jq^{\binom{j}{2}} &= \sum_{d=0}^c q^{d(a-b-c+d)} \qbinom{a-b}{c-d} \sum_{j=0}^b \qbinom{b}{j} \qbinom{b-j}{d} (-1)^jq^{\binom{j}{2}} \\
	&= \sum_{d=0}^c q^{d(a-b-c+d)} \qbinom{a-b}{c-d} \qbinom{b}{d} \sum_{j=0}^b \qbinom{b-d}{j} (-1)^jq^{\binom{j}{2}} \\
	&= \sum_{d=0}^c q^{d(a-b-c+d)} \qbinom{a-b}{c-d} \qbinom{b}{d} (1;q)_{b-d} \\
	&= q^{b(a-c)} \qbinom{a-b}{c-b},
\end{align*}
as desired.
\end{proof}

For $0\leqslant k,\ell\leqslant n$, we have
\begin{align}
	F_{k,\ell} &= \sum_{j,m=0}^n C_{k,j}G_{j,m}(C^{-1})_{m,\ell} \label{F expression} \\
	&= \sum_{m=0}^n \qbinom{m}{\ell}(-1)^{m-\ell}\sigma^m q^{\binom{m}{2}+\binom{m-\ell}{2}} \sum_{j=0}^n \qbinom{k}{j}\qbinom{n-j}{n-m} (-1)^j q^{\binom{j}{2}} \notag \\
	&= \sum_{m=0}^n \qbinom{n-k}{m}\qbinom{m}{\ell}(-1)^{m-\ell}\sigma^m q^{km+\binom{m}{2}+\binom{m-\ell}{2}} \notag \\
	&= \qbinom{n-k}{\ell} \sigma^{\ell}q^{k\ell+\binom{\ell}{2}} \sum_{h=0}^{n-k-\ell}\qbinom{n-k-\ell}{h} (-1)^h \sigma^h q^{h(h+k+\ell-1)} \notag
\end{align}
by Lemma \ref{technical identity} (with $(a,b,c)=(n,k,n-m)$), where we set $h=m-\ell$ in the last line above.
In particular, it follows that $F$ is upper anti-triangular, i.e., $F_{k,\ell}=0$ whenever $k+\ell>n$.
Moreover, it is immediate to see that
\begin{equation*}
	F_{k,\ell}\cdot\phi_k\qbinom{n}{k}=F_{\ell,k}\cdot\phi_{\ell}\qbinom{n}{\ell} \qquad (0\leqslant k,\ell\leqslant n).
\end{equation*}
Thus, if we define the matrix $A'$ in \eqref{A'} by \eqref{our choice for A'}, then $A'$ is symmetric since
\begin{equation*}
	a_{k,\ell}'=\frac{F_{\ell,k}}{1+\sigma}\cdot\frac{\phi_{\ell}\qbinom{n}{\ell}}{\phi_k\qbinom{n}{k}}\cdot\frac{\phi_k^{\frac{1}{2}}\qbinom{n}{k}^{\frac{1}{2}}}{\phi_{\ell}^{\frac{1}{2}}\qbinom{n}{\ell}^{\frac{1}{2}}}=a_{\ell,k}' \qquad (0\leqslant k\leqslant n,\, 0\leqslant\ell\leqslant n-k).
\end{equation*}

It seems that the entries of $F$ have no simpler expression in general.
However, if we let $q\rightarrow 1$, then
\begin{equation*}
	\lim_{q\rightarrow 1} G_{k,\ell}=\binom{n-k}{n-\ell}(-1)^k\sigma^{\ell} \qquad (0\leqslant k\leqslant\ell\leqslant n),
\end{equation*}
and we also have
\begin{equation*}
	\lim_{q\rightarrow 1} F_{k,\ell} = \binom{n-k}{\ell} \sigma^{\ell}(1-\sigma)^{n-k-\ell} \qquad (0\leqslant k\leqslant n,\, 0\leqslant\ell\leqslant n-k).
\end{equation*}
The relation \eqref{G and F} after taking the limit $q\rightarrow 1$ was shown earlier in \cite[Lemmas 2.9, 2.21]{OSST2014pre}.
We also note that the above limit of $F$ is closely related to the matrix $A^{(n)}$ (with $p=\frac{\sigma}{1+\sigma}$) considered in \cite{Friedgut2008C}.

A $1$-eigenvector of $G$ is given by $(1,0,\dots,0)^{\mathsf{T}}\in\mathbb{R}^{n+1}$, since $G$ is upper triangular and $G_{0,0}=1$.
Then, a $1$-eigenvector of $F=CGC^{-1}$ is
\begin{equation*}
	C(1,0,\dots,0)^{\mathsf{T}}=(1,1,\dots,1)^{\mathsf{T}}.
\end{equation*}

It follows from \eqref{block-diagonalization} that the matrix $A'$ in \eqref{A'} satisfies
\begin{equation*}
	A'\bm{u}_{0,1}^{\ell} = \sum_{k=0}^{n-\ell} a_{k,\ell}' \bm{u}_{0,1}^k \qquad (0\leqslant \ell\leqslant n).
\end{equation*}
In other words, the matrix $A_0'\in\mathbb{R}^{(n+1)\times(n+1)}$ representing the action of $A'$ on the subspace $\bm{V}_{0,1}$ with respect to the basis $\bm{u}_{0,1}^0,\bm{u}_{0,1}^1,\dots,\bm{u}_{0,1}^n$ is upper anti-triangular and is given by
\begin{equation*}
	(A_0')_{k,\ell}=\begin{cases} a_{k,\ell}' & \text{if} \ k+\ell\leqslant n, \\ 0 & \text{if} \ k+\ell>n, \end{cases} \qquad (0\leqslant k,\ell\leqslant n).
\end{equation*}
From now on, we define $A'$ by \eqref{our choice for A'}.
Then, since $F$ is also upper anti-triangular, we have
\begin{equation}\label{A0'}
	A_0' = \frac{1}{1+\sigma} D_0FD_0^{-1},
\end{equation}
where $D_0\in\mathbb{R}^{(n+1)\times(n+1)}$ is the diagonal matrix with diagonal entries $(D_0)_{k,k}=\phi_k^{\frac{1}{2}}\qbinom{n}{k}^{\frac{1}{2}}$ $(0\leqslant k\leqslant n)$.
In particular, $A_0'$ has eigenvalues
\begin{equation*}
	\frac{(-1)^k\sigma^kq^{k(k-1)}}{1+\sigma} \qquad (0\leqslant k\leqslant n).
\end{equation*}
Moreover, it follows from the above comment that a $\frac{1}{1+\sigma}$-eigenvector of $A_0'$ is given by
\begin{equation}\label{w0}
	\bm{w}_0:=D_0(1,1,\dots,1)^{\mathsf{T}}=\left(\phi_0^{\frac{1}{2}}\qbinom{n}{0}^{\frac{1}{2}}, \phi_1^{\frac{1}{2}}\qbinom{n}{1}^{\frac{1}{2}},\dots, \phi_n^{\frac{1}{2}}\qbinom{n}{n}^{\frac{1}{2}}\right)^{\!\!\mathsf{T}}.
\end{equation}
By \eqref{special eigenvector}, the matrix representing the action of $\Delta^{\frac{1}{2}}J\Delta^{\frac{1}{2}}=(\Delta^{\frac{1}{2}}\bm{1})(\Delta^{\frac{1}{2}}\bm{1})^{\mathsf{T}}$ on the subspace $\bm{V}_{0,1}$ with respect to the same basis $\bm{u}_{0,1}^0,\bm{u}_{0,1}^1,\dots,\bm{u}_{0,1}^n$ is $\bm{w}_0(\bm{w}_0)^{\mathsf{T}}$.
Indeed, we have
\begin{equation*}
	\Delta^{\frac{1}{2}}J\Delta^{\frac{1}{2}} \bm{u}_{0,1}^{\ell} = \phi_{\ell}^{\frac{1}{2}} \qbinom{n}{\ell}^{\frac{1}{2}} \Delta^{\frac{1}{2}}\bm{1} = \sum_{k=0}^n \phi_k^{\frac{1}{2}} \qbinom{n}{k}^{\frac{1}{2}}\phi_{\ell}^{\frac{1}{2}} \qbinom{n}{\ell}^{\frac{1}{2}} \bm{u}_{0,1}^k \qquad (0\leqslant\ell\leqslant n).
\end{equation*}
Since
\begin{equation*}
	(\bm{w}_0)^{\mathsf{T}}\bm{w}_0=\sum_{k=0}^n \phi_k\qbinom{n}{k}=\mu_{\sigma}(\Omega_n)=1,
\end{equation*}
the matrix $\bm{w}_0(\bm{w}_0)^{\mathsf{T}}$ has $\bm{w}_0$ as a $1$-eigenvector.
Since $\Delta^{\frac{1}{2}}J\Delta^{\frac{1}{2}}$ is a rank-one matrix, this is the only nontrivial action of $\Delta^{\frac{1}{2}}J\Delta^{\frac{1}{2}}$, i.e., all the other eigenvalues are zero.
Recall our choice of $\alpha$ and $Z'$.
The vector $\bm{w}_0$ is an eigenvector of the action of
\begin{equation*}
	S'=\frac{\sigma}{1+\sigma}I-\Delta^{\frac{1}{2}}J\Delta^{\frac{1}{2}}+A'
\end{equation*}
on $\bm{V}_{0,1}$ with eigenvalue
\begin{equation*}
	\frac{\sigma}{1+\sigma}-1+\frac{1}{1+\sigma}=0.
\end{equation*}
The other $n$ eigenvalues of $S'$ on $\bm{V}_{0,1}$ are given by
\begin{equation}\label{0th block eigenvalues}
	\frac{\sigma}{1+\sigma}+\frac{(-1)^k\sigma^kq^{k(k-1)}}{1+\sigma} \qquad (1\leqslant k\leqslant n).
\end{equation}

Next, we consider the actions of $S'$ on the other subspaces $\bm{V}_{i,r}$, where $1\leqslant i\leqslant \lfloor \frac{n}{2}\rfloor$ and $1\leqslant r\leqslant d_i$.
By \eqref{block-diagonalization}, the matrix $A_i'$, indexed by $i,i+1,\dots,n-i$, representing the action of $A'$ on $\bm{V}_{i,r}$ with respect to the basis $\bm{u}_{i,r}^i,\bm{u}_{i,r}^{i+1},\dots,\bm{u}_{i,r}^{n-i}$ is given by
\begin{equation*}
	(A_i')_{k,\ell} = a_{k,\ell}'\cdot\frac{\theta_i^{k,\ell}}{\theta_0^{k,\ell}} = a_{k,\ell}' \cdot (-1)^i q^{\binom{i}{2}-\frac{i(k+\ell)}{2}} \!\left( \frac{ (q;q)_{n-k-i} (q;q)_{n-\ell-i} (q;q)_k (q;q)_{\ell} }{ (q;q)_{n-k} (q;q)_{n-\ell} (q;q)_{k-i} (q;q)_{\ell-i} } \right)^{\!\frac{1}{2}}
\end{equation*}
for $i\leqslant k,\ell\leqslant n-i$.
We then have (cf.~\eqref{our choice for A'})
\begin{equation}\label{Ai'}
	A_i'=\frac{1}{1+\sigma} D_iF_iD_i^{-1},
\end{equation}
where
\begin{equation*}
	(F_i)_{k,\ell}=F_{k,\ell} \cdot (-1)^i q^{\binom{i}{2}-ik} \frac{ (q;q)_{n-k-i} (q;q)_{\ell} }{ (q;q)_{n-k} (q;q)_{\ell-i} } \qquad (i\leqslant k,\ell\leqslant n-i),
\end{equation*}
and $D_i$ is diagonal with diagonal entries
\begin{equation*}
	(D_i)_{k,k}= \phi_k^{\frac{1}{2}}\qbinom{n}{k}^{\frac{1}{2}} \cdot q^{\frac{ik}{2}} \!\left(\frac{(q;q)_{n-k}(q;q)_k}{(q;q)_{n-k-i}(q;q)_{k-i}} \right)^{\!\frac{1}{2}} \qquad (i\leqslant k\leqslant n-i).
\end{equation*}
If we write $F=F_{n;\,\sigma}$ to specify the parameters, then using \eqref{F expression},
\begin{equation*}
	\qbinom{n-k}{\ell}=\qbinom{n-k-i}{\ell-i}\frac{(q;q)_{n-k}(q;q)_{\ell-i}}{(q;q)_{n-k-i}(q;q)_{\ell}},
\end{equation*}
and
\begin{equation*}
	k\ell=(k-i)(\ell-i)+i(k+\ell)-i^2, \qquad \binom{\ell}{2}=\binom{\ell-i}{2}+i\ell-\binom{i+1}{2},
\end{equation*}
it is routinely verified that
\begin{equation}\label{different parameters}
	(F_i)_{k,\ell}=(F_{n-2i;\,\sigma q^{2i}})_{k-i,\ell-i} \cdot (-1)^i \sigma^i q^{i(i-1)} \qquad (i\leqslant k,\ell\leqslant n-i),
\end{equation}
where we note that the rows and columns of $F_{n-2i;\,\sigma q^{2i}}$ are indexed by $0,1,\dots,n-2i$.
It follows that the eigenvalues of $A'$ on $\bm{V}_{i,r}$ are given by
\begin{equation*}
	(-1)^{h}(\sigma q^{2i})^hq^{h(h-1)} \cdot \frac{(-1)^i \sigma^i q^{i(i-1)}}{1+\sigma} = \frac{(-1)^k\sigma^kq^{k(k-1)}}{1+\sigma} \qquad (i\leqslant k\leqslant n-i)
\end{equation*}
where $h=k-i$, and therefore those of $S'$ are
\begin{equation}\label{ith block eigenvalues}
	\frac{\sigma}{1+\sigma}+\frac{(-1)^k\sigma^kq^{k(k-1)}}{1+\sigma} \qquad (i\leqslant k\leqslant n-i).
\end{equation}

For the matrix $S'$ to be positive semidefinite, all the eigenvalues in \eqref{0th block eigenvalues} and \eqref{ith block eigenvalues} must be nonnegative.
This is equivalent to
\begin{equation*}
	\sigma^kq^{k(k-1)} \leqslant \sigma \qquad (1\leqslant k\leqslant n,\ k:\text{odd}),
\end{equation*}
which then simplifies to the condition given in Theorem \ref{thm1} (when $n\geqslant 3$).
If this condition is satisfied, then the matrix $S=\Delta^{\frac{1}{2}} S'\Delta^{\frac{1}{2}}$ gives a feasible solution to (D) with objective value $\frac{\sigma}{1+\sigma}$, which is attained by $A_n^{(1)}$ defined by \eqref{def:An}.

For the rest of the proof, assume that $\sigma< q^{-2\lfloor \frac{n-1}{2}\rfloor-1}$.
Let $U\subset\Omega_n$ be an intersecting family such that $\mu_{\sigma}(U)=\frac{\sigma}{1+\sigma}$.
Let $\bm{x}\in\mathbb{R}^{\Omega_n}$ be the characteristic vector of $U$, and let the matrix $X\in S\mathbb{R}^{\Omega_n\times\Omega_n}$ be as in \eqref{X from U}.
Then, equality is attained in \eqref{weak duality}, and hence it follows that $\operatorname{tr}(SX)=0$, or equivalently, $S'\Delta^{\frac{1}{2}}\bm{x}=0$.

Recall that $\bm{w}_0$ is a $0$-eigenvector of the action of $S'$ on $\bm{V}_{0,1}$.
The corresponding $0$-eigenvector of $S'$ (in $\mathbb{R}^{\Omega_n}$) is (cf.~\eqref{special eigenvector})
\begin{equation}\label{v0}
	\bm{v}_0:=\Delta^{\frac{1}{2}}\bm{1}=\sum_{k=0}^n \phi_k^{\frac{1}{2}}\qbinom{n}{k}^{\frac{1}{2}} \bm{u}_{0,1}^k.
\end{equation}
The eigenvalues of $S'$ in \eqref{ith block eigenvalues} are zero if and only if $(i,k)=(1,1)$, in which case, we can similarly see that the corresponding $0$-eigenvectors of $S'$ are of the form
\begin{equation*}
	\bm{v}_r=\sum_{k=1}^{n-1} \eta_k \bm{u}_{1,r}^k \qquad (1\leqslant r\leqslant d_1),
\end{equation*}
where $\eta_k\ne 0$ for all $k$.
More specifically,
\begin{equation*}
	(\eta_1,\eta_2,\dots,\eta_{n-1})^{\mathsf{T}}=D_1(1,1,\dots,1)^{\mathsf{T}}=\big((D_1)_{1,1},(D_1)_{2,2},\dots,(D_1)_{n-1,n-1}\big)^{\!\mathsf{T}}.
\end{equation*}
See \eqref{Ai'} and \eqref{different parameters}.
On the other hand, the eigenvalues of $S'$ in \eqref{0th block eigenvalues} are zero if and only if $k=1$.
Since $G$ is upper triangular, and $G_{0,0}=1$ and $G_{1,1}=-\sigma$, a $(-\sigma)$-eigenvector of $G$ is given by $(\nu,1,0,\dots,0)^{\mathsf{T}}\in\mathbb{R}^{n+1}$, where $\nu=-\frac{G_{0,1}}{1+\sigma}$.
Then, a $(-\sigma)$-eigenvector of $F=CGC^{-1}$ is
\begin{equation*}
	C(\nu,1,0,\dots,0)^{\mathsf{T}}=\left(\nu+\qbinom{0}{1},\nu+\qbinom{1}{1},\nu+\qbinom{2}{1},\dots,\nu+\qbinom{n}{1}\right)^{\!\!\mathsf{T}},
\end{equation*}
and the corresponding $0$-eigenvector of $S'$ becomes (cf.~\eqref{A0'})
\begin{equation*}
	\bm{v}_0'=\sum_{k=0}^n \phi_k^{\frac{1}{2}}\qbinom{n}{k}^{\frac{1}{2}}\!\left(\nu+\qbinom{k}{1}\right) \bm{u}_{0,1}^k.
\end{equation*}
Since $S'\Delta^{\frac{1}{2}}\bm{x}=0$, the vector $\Delta^{\frac{1}{2}}\bm{x}$ must be a linear combination of the above vectors:
\begin{equation*}
	\Delta^{\frac{1}{2}}\bm{x}= c_0\bm{v}_0+c_1\bm{v}_1+\cdots+c_{d_1}\bm{v}_{d_1}+c_0'\bm{v}_0'.
\end{equation*}
Expand the RHS above in terms of the $\bm{u}_{i,r}^k$.
Note that $U$ does not contain $0\in\ksubsp{n}{0}$, so that the coefficient of $\bm{u}_{0,1}^0$ is zero, i.e., $c_0+c_0'\nu=0$.
Suppose now that $U\cap\ksubsp{n}{1}=\emptyset$.
Then, the coefficients of $\bm{u}_{0,1}^1$ and $\bm{u}_{1,r}^1$ $(1\leqslant r\leqslant d_1)$ are all zero because these vectors form a basis of $\mathbb{R}^{\ksubsp{n}{1}}$.
That the coefficient of $\bm{u}_{0,1}^1$ equals zero is equivalent to $c_0+c_0'(\nu+1)=0$.
Combining this with $c_0+c_0'\nu=0$, we have $c_0=c_0'=0$.
Moreover, the coefficient of $\bm{u}_{1,r}^1$ equals $c_r\eta_1$ for $1\leqslant r\leqslant d_1$, and hence $c_1=\cdots=c_{d_1}=0$ since $\eta_1\ne 0$.
It follows that $\Delta^{\frac{1}{2}}\bm{x}=0$, which is absurd.
We have now shown that $U\cap\ksubsp{n}{1}\ne\emptyset$.
Since $U$ is a maximal intersecting family, we must have $U=A_n^{(1)}$ for some $y\in\ksubsp{n}{1}$.
This completes the proof of Theorem \ref{thm1}.

\section{Proof of Theorem~\ref{thm2}}\label{sec3}

In this section, we abbreviate
\begin{equation*}
	\mu:=\mu_{\theta,n}, \qquad \sigma:=\sigma_{\theta,n},
\end{equation*}
except in the statement of lemmas.
We will also write $\phi=\phi_{\theta,n}:=\phi_{\sigma_{\theta,n},n}$.
\begin{lemma}\label{la2}
If $0<\theta<\frac{1}{2}$, then
\begin{equation*}
	\lim_{n\to\infty}\mu_{\theta,n}(A_n^{(t)})^{\frac{1}{n}}=q^{-(1-\theta)t}.
\end{equation*}
\end{lemma}
\begin{proof}
By \eqref{mu_s(A)}, we have
\begin{equation}\label{mu(A)}
	\mu(A_n^{(t)})=\prod_{j=0}^{t-1}\left(1+q^{(1-\theta)n-j}\right)^{-1}.
\end{equation}
Since $(1+q^{(1-\theta)n-j})^{-\frac{1}{n}}\to q^{-(1-\theta)}$ for $0\leqslant j\leqslant t-1$, we have
$\mu(A_n^{(t)})^{\frac{1}{n}}\to q^{-(1-\theta)t}$.
\end{proof}

The next result shows that the distribution 
\begin{equation*}
	\Phi(k)=\Phi_{\theta,n}(k):=\qbinom{n}{k}\phi_{\theta,n}(k) \qquad (k=0,1,\ldots,n)
\end{equation*}
concentrates around $k\sim n\theta$.

\begin{lemma}\label{la3}
Let $0<\theta<1$.
For every $\epsilon>0$, there exists $L>1$ such that  
\begin{equation*}
	\sum_{|k-\theta n|>L}\!\!\Phi_{\theta,n}(k)<\epsilon
\end{equation*}
for sufficiently large $n$, where the sum is over all $k=0,1,\dots,n$ such that $|k-\theta n|>L$.
\end{lemma}
\begin{proof}
Define a probability measure $\Psi=\Psi_{\theta,n}$ on the points 
$q^{k-\theta n}$ $(k = 0,1,\ldots,n)$ by
\begin{equation*}
	\Psi(q^{k-\theta n})=\Phi(k)\qquad(k=0,1,\ldots,n).
\end{equation*}
By \eqref{q-binom}, the mean is computed as
\begin{equation*}
	\mathbb{E}[X]= \sum_{k=0}^n q^{k-\theta n}\, \Phi(k)
	=\frac{(-\sigma q;q)_n}{q^{\theta n}(-\sigma;q)_n}
	=\frac{1+\sigma q^n}{q^{\theta n}(1+\sigma)},
\end{equation*}
which converges to $1$ when $n\to\infty$.
Also, 
\begin{equation*}
	\mathbb{E}[X^2]=
	\sum_{k=0}^n q^{2k-2\theta n}\, \Phi(k)
	=\frac{(-\sigma q^2;q)_n}{q^{2\theta n}(-\sigma;q)_n}
	=\frac{(1+\sigma q^n)(1+\sigma q^{n+1})}{q^{2\theta n}(1+\sigma)(1+\sigma q)},
\end{equation*}
which converges to $q$ when $n\to\infty$.
Hence, the variance satisfies
\begin{equation*}
	\mathbb{V}[X]=\mathbb{E}[X^2]-\mathbb{E}[X]^2 \to q-1.
\end{equation*}

In what follows, let $n$ be sufficiently large so that 
$\mathbb{E}[X] < 2$ and $\mathbb{V}[X] < q$.
There exists $L > 1$ such that $(q^L-2)^2\epsilon>2q$.
Then, we have
\begin{equation*}
	q>\mathbb{V}[X]\geqslant\sum_{k>\theta n+L}\!(q^{k-\theta n}-\mathbb{E}[X])^2\,\Phi(k)
	\geqslant(q^L-2)^2\sum_{k>\theta n+L}\!\!\Phi(k),
\end{equation*}
from which it follows that
\begin{equation*}
	\sum_{k>\theta n+L}\!\!\Phi(k)<\frac{q}{(q^L-2)^2}<\frac{\epsilon}{2}.
\end{equation*}

Next, we have
\begin{equation*}
	\mathbb{E}[X^{-1}] = \frac{q^{\theta n}(1+\sigma q^{-1})}{1+\sigma q^{n-1}}\to q,
	\qquad
	\mathbb{E}[X^{-2}] = \frac{q^{2\theta n}(1+\sigma q^{-2})(1+\sigma q^{-1})}
	{(1+\sigma q^{n-2})(1+\sigma q^{n-1})}\to q^3,
\end{equation*}
so
\begin{equation*}
	\mathbb{V}[X^{-1}] = \mathbb{E}[X^{-2}] - \mathbb{E}[X^{-1}]^2\to q^3 - q^2.
\end{equation*}
By a similar argument, we can show that there exists $L' > 1$ such that
\begin{equation*}
	\sum_{k<\theta n-L'}\!\!\Phi(k)<\frac{\epsilon}{2},
\end{equation*}
for sufficiently large $n$.
The result now follows by replacing $L$ by $L'$ if $L'>L$.
\end{proof}

For our purpose, we need a stronger tail bound.

\begin{claim}\label{cl4}
Let $0<\theta<\frac{1}{2}$.
We have
\begin{equation}\label{claim 4 evaluation}
	\sum_{k>\frac{n}{2}}\Phi_{\theta,n}(k)=o(q^{-(1-\theta)tn}),
\end{equation}
where the sum is over all integers $k$ with $\frac{n}{2}<k\leqslant n$.
\end{claim}

\begin{proof}
We write $m:=\lceil\frac{n}{2}\rceil$ and $s:=\lceil\theta n\rceil$ for typographical reasons.
First, we claim that
\begin{equation}\label{phi(n/2)}
	\Phi(m)<q^{-\frac{1}{2}(m-s)^2+O(n)}=q^{-\frac{1}{2}(\frac{1}{2}-\theta)^2n^2+O(n)}.
\end{equation}
In view of Lemma~\ref{la3}, we estimate (cf.~\eqref{q-binom})
\begin{equation}\label{eq1b}
	\Phi(m)
	=\frac{\qbinom{n}{m} \sigma^m q^{\binom{m}{2}}}
	{\sum_{k=0}^n\qbinom{n}{k}\sigma^{k}q^{\binom{k}{2}}}
	<\frac{\qbinom{n}{m} \sigma^m q^{\binom{m}{2}}}
	{\qbinom{n}{s}\sigma^sq^{\binom{s}{2}}}
	=\frac{\qbinom{n}{m}}{\qbinom{n}{s}}
	\sigma^{m-s} q^{\binom{m}{2} - \binom{s}{2}}.
\end{equation}
For further estimation of the RHS above, we note that
\begin{equation*}
	\frac{[n-m]}{[s]}=\frac{q^{n-m}-1}{q^s-1}<q^{n-m-s+1}=q^{m-s+O(1)},
\end{equation*}
and so
\begin{equation}\label{eq2b}
	\frac{\qbinom{n}{m}}{\qbinom{n}{s}}
	=\frac{[s]![n-s]!}{[m]![n-m]!}
	=\prod_{j=1}^{m-s}\frac{[n-m+j]}{[s+j]}
	<\left(\frac{[n-m]}{[s]}\right)^{\!\! m-s}
	\!\! <q^{(m-s)^2+O(n)}
\end{equation}
for sufficiently large $n$.
We also have
\begin{equation}\label{eq3b}
	\sigma^{m-s} q^{\binom{m}{2} - \binom{s}{2}}
	=\left(q^{-(1-\theta)n}\right)^{m-s}q^{\frac{m^2}{2}-\frac{s^2}{2}+O(n)}
	=q^{-\frac{3}{2}(m-s)^2+O(n)}.
\end{equation}
Substituting \eqref{eq2b} and \eqref{eq3b} into the RHS of \eqref{eq1b}, 
we get \eqref{phi(n/2)}.

Next, we verify that $\Phi(k)$ is decreasing in $k$ for $k\geqslant m$.
Indeed, since $\frac{[n-k]}{[k+1]}<q^{n-2k}$ as above, it follows that
\begin{equation*}
	\frac{\Phi(k+1)}{\Phi(k)}=\frac{\qbinom{n}{k+1}\sigma^{k+1}q^{\binom{k+1}{2}}}
	{\qbinom{n}{k}\sigma^{k}q^{\binom{k}{2}}}=\frac{[n-k]}{[k+1]}\sigma q^k
	<q^{\theta n-k}\leqslant q^{\theta n-m}<1.
\end{equation*}
Hence, it follows from \eqref{phi(n/2)} that 
\begin{equation*}
	\sum_{k>\frac{n}{2}}\Phi(k)
	<m\Phi(m)<q^{-\frac{1}{2}(\frac{1}{2}-\theta)^2n^2+O(n)}
	=o(q^{-(1-\theta)tn})
\end{equation*}
for sufficiently large $n$.
\end{proof}

\begin{rem}\label{linear}
The RHS of \eqref{claim 4 evaluation} can be replaced by $o(q^{-Rn})$ for any fixed $R>0$.
For our purpose, we need $R\geqslant (1-\theta)t$.
\end{rem}

We now invoke the following result.

\begin{oldthm}[Frankl--Wilson \cite{FW1986}]\label{thmD}
Let $n\geqslant 2k$.
If a family $U\subset\ksubsp{n}{k}$ is $t$-intersecting, then
$|U|\leqslant\qbinom{n-t}{k-t}$.
\end{oldthm}

For the characterization of the optimal families in Theorem \ref{thmD}, see \cite{Hsieh1975DM,Tanaka2006JCTA}.

\begin{lemma}\label{la5}
 If $0<\theta<\frac{1}{2}$, then 
$\lim_{n\to\infty}f(n,t,\sigma_{\theta,n})=q^{-(1-\theta)t}$.
\end{lemma}

\begin{proof}
By Lemma~\ref{la2}, we have 
\begin{equation*}
	f(n,t,\sigma)\geqslant \mu(A_n^{(t)})^{\frac{1}{n}}\rightarrow q^{-(1-\theta)t}.
\end{equation*}
On the other hand, let $U_n\subset\Omega_n$ be a $t$-intersecting family satisfying
$\mu(U_n)=f(n,t,\sigma)^n$, and let $U_n^{(k)}=U_n\cap\ksubsp{n}{k}$.
We have
\begin{equation*}
	\mu(U_n)=\sum_{k=0}^n|U_n^{(k)}|\phi(k)
\leqslant\sum_{k\leqslant\frac{n}{2}}|U_n^{(k)}|\phi(k)+\sum_{k>\frac{n}{2}}\Phi(k).
\end{equation*}
By Theorem~\ref{thmD} and \eqref{mu(A)}, 
\begin{equation*}
	\sum_{k\leqslant\frac{n}{2}}|U_n^{(k)}|\phi(k)\leqslant\sum_{k\leqslant\frac{n}{2}}\qbinom{n-t}{k-t}\phi(k) \leqslant \mu(A_n^{(t)})=\prod_{j=0}^{t-1}\left(1+q^{(1-\theta)n-j}\right)^{-1},
\end{equation*}
and by Claim~\ref{cl4},
\begin{equation*}
	\sum_{k>\frac{n}{2}}\Phi(k) = o(q^{-(1-\theta)tn}).
\end{equation*}
Thus, we have
\begin{align*}
	\mu(U_n) &\leqslant \mu(A_n^{(t)})+o(q^{-(1-\theta)tn}) \\
	&=\mu(A_n^{(t)})\!\left(1+o\!\left(\prod_{j=0}^{t-1}(q^{-(1-\theta)n}+q^{-j})\right)\!\!\right) \\
	&=\mu(A_n^{(t)})(1+o(1)).
\end{align*}
Finally, it follows from Lemma~\ref{la2} that
\begin{equation*}
	f(n,t,\sigma)= \mu(U_n)^{\frac{1}{n}}\leqslant\mu(A_n^{(t)})^{\frac{1}{n}}
	(1+o(1))^{\frac{1}{n}}\to q^{-(1-\theta)t}. \qedhere
\end{equation*}
\end{proof}

\begin{claim}\label{cl6}
Let $\frac{1}{2}<\theta<1$.
There exists $\delta>0$ such that
\begin{equation*}
	\sum_{k<\frac{n+t}{2}} \!\Phi_{\theta,n}(k)<q^{-\delta n^2}
\end{equation*}
for sufficiently large $n$, where the sum is over all integers $k$ with $0\leqslant k<\frac{n+t}{2}$.
\end{claim}
\begin{proof}
The proof is similar to that of Claim~\ref{cl4}.
We write $m':=\lceil\frac{n+t}{2}\rceil$ and $s:=\lceil\theta n\rceil$.
We first claim that
\begin{equation}\label{phi(n+t/2)}
	\Phi(m')<q^{-\frac{1}{2}(s-m')^2+O(n)}=q^{-\frac{1}{2}(\theta-\frac{1}{2})^2n^2+O(n)}.
\end{equation}
We have (cf.~\eqref{q-binom})
\begin{equation*}
	\Phi(m')=\frac{\qbinom{n}{m'}\sigma^{m'}q^{\binom{m'}{2}}}{\sum_{k=0}^n\qbinom{n}{k}\sigma^{k}q^{\binom{k}{2}}}
	<\frac{\qbinom{n}{m'} \sigma^{m'} q^{\binom{m'}{2}}}
	{\qbinom{n}{s}\sigma^s q^{\binom{s}{2}}}
	=\frac{\qbinom{n}{m'}}{\qbinom{n}{s}}
	\sigma^{m'-s} q^{\binom{m'}{2} - \binom{s}{2}}.
\end{equation*}
Since $\frac{[m']}{[n-s]}<q^{m'-n+s+1}=q^{s-m'+O(1)}$, it follows that
\begin{equation*}
	\frac{\qbinom{n}{m'}}{\qbinom{n}{s}}
	=\frac{[s]![n-s]!}{[m']![n-m']!}
	=\prod_{j=1}^{s-m'} \frac{[m'+j]}{[n-s+j]}
	<\left( \frac{[m']}{[n-s]}\right)^{\!\! s-m'}
	\!\! =q^{(s-m')^2+O(n)}
\end{equation*}
for sufficiently large $n$.
We also have 
\begin{equation*}
	\sigma^{m'-s} q^{\binom{m'}{2} - \binom{s}{2}}
	=\left(q^{-(1-\theta)n}\right)^{\! m'-s}q^{\frac{(m')^2}{2}-\frac{s^2}{2}+O(n)}=q^{-\frac{3}{2}(s-m')^2+O(n)}.
\end{equation*}
Hence, we get \eqref{phi(n+t/2)}.

Next, we verify that $\Phi(k)$ is increasing in $k$ for $k\leqslant m'$, provided that $n$ is sufficiently large.
This follows from 
$\frac{[k]}{[n-k+1]}<q^{2k-n}$ and
\begin{equation*}
	\frac{\Phi(k-1)}{\Phi(k)}=\frac{\qbinom{n}{k-1}\sigma^{k-1}q^{\binom{k-1}{2}}}
	{\qbinom{n}{k}\sigma^{k}q^{\binom{k}{2}}}=\frac{[k]}{[n-k+1]\sigma q^{k-1}}
	<q^{k-\theta n+1}\leqslant q^{m'-\theta n+1}<1.
\end{equation*}
Hence, if we choose $\delta$ such that $0<\delta<\frac{1}{2}(\theta-\frac{1}{2})^2$, then we have
\begin{equation*}
	\sum_{k<\frac{n+t}{2}}\!\Phi(k)<m'\Phi(m')<q^{-\delta n^2}
\end{equation*}
for sufficiently large $n$.
\end{proof}

\begin{lemma}\label{la7}
 If $\frac{1}{2}<\theta<1$, then 
$\lim_{n\to\infty}f(n,t,\sigma_{\theta,n})=1$.
\end{lemma}
\begin{proof}
Clearly, we have $f(n,t,\sigma)\leqslant 1$.
So we need to show that $\varliminf_{n\to\infty}f(n,t,\sigma)\geqslant 1$.
Define a $t$-intersecting family $B_n$ by
\begin{equation*}
	B_n:=\left\{x\in\Omega_n:\dim x\geqslant\frac{n+t}{2}\right\}.
\end{equation*}
By Claim~\ref{cl6}, we have 
\begin{equation*}
	\mu(B_n)= 1-\sum_{k<\frac{n+t}{2}}\Phi(k)>1-q^{-\delta n^2}
\end{equation*}
for sufficiently large $n$.
Then, the desired inequality follows from
\begin{equation*}
	f(n,t,\sigma)\geqslant\mu(B_n)^{\frac{1}{n}}\geqslant
	\left(1-q^{-\delta n^2}\right)^{\frac{1}{n}}\rightarrow 1. \qedhere
\end{equation*}
\end{proof}

\begin{proof}[Proof of Theorem~\ref{thm2}]
Immediate from Lemma~\ref{la5} and Lemma~\ref{la7}.
\end{proof}

\section{Proof of Theorem~\ref{thm3}}\label{sec4}

Here, we list some results concerning cross $t$-intersecting families
of uniform subspaces.
We omit the descriptions of the optimal families.

\begin{oldthm}[Tokushige \cite{T2013}]\label{thmE}
Let $n\geqslant 2k$.
If $U\subset\ksubsp{n}{k}$ and $W\subset\ksubsp{n}{k}$ are cross $t$-intersecting, then
$|U||W|\leqslant\qbinom{n-t}{k-t}^2$.
\end{oldthm}

\begin{oldthm}[Suda--Tanaka \cite{ST2014BLMS}]\label{thmF}
Let $n\geqslant 2k$ and $n\geqslant 2\ell$.
If $U\subset\ksubsp{n}{k}$ and
$W\subset\ksubsp{n}{\ell}$ are cross $1$-intersecting, then
$|U||W|\leqslant\qbinom{n-1}{k-1}\qbinom{n-1}{\ell-1}$.
\end{oldthm}

\begin{oldthm}[Cao--Lu--Lv--Wang \cite{CLLW}]\label{thmG}
Let $n\geqslant k+\ell+t+1$.
If $U\subset\ksubsp{n}{k}$ and
$W\subset\ksubsp{n}{\ell}$ are cross $t$-intersecting, then
$|U||W|\leqslant\qbinom{n-t}{k-t}\qbinom{n-t}{\ell-t}$.
\end{oldthm}

For the proof of Theorem~\ref{thm3}, we use Theorem~\ref{thmG}.
Note that while Theorem~\ref{thmG} is the most general result so far, 
it does not fully contain Theorem~\ref{thmE} and Theorem~\ref{thmF}.

The next claim can be shown exactly in the same way (with slightly more
cumbersome computation) as Claim~\ref{cl4}, and we omit the proof.
See also Remark \ref{linear}.

\begin{claim}\label{cl8}
Let $0<\theta<\frac{1}{2}$, and let $t$ be a fixed positive integer.
Then, we have
\begin{equation*}
	\sum_{k>\frac{n-t-1}{2}}\!\!\Phi_{\theta,n}(k)= o(q^{-2tn}),
\end{equation*}
where the sum is over all integers $k$ with $\frac{n-t-1}{2}<k\leqslant n$.
\end{claim}

\begin{proof}[Proof of Theorem~\ref{thm3}]
For $i=1,2$, we write $\sigma_i:=\sigma_{\theta_i,n}$, $\phi_i:=\phi_{\theta_i,n}$, and $\mu_i:=\mu_{\theta_i,n}$ for brevity.
Suppose that cross $t$-intersecting families $U_n,W_n\subset\Omega_n$
satisfy
\begin{equation*}
	g(n,t,\sigma_1,\sigma_2)^n=\mu_1(U_n)\mu_2(W_n),
\end{equation*}
and let $U_n^{(k)}=U_n\cap\ksubsp{n}{k}$ and $W_n^{(\ell)}=W_n\cap\ksubsp{n}{\ell}$.
If $k,\ell\leqslant\frac{n-t-1}{2}$, then $n\geqslant k+\ell+t+1$, 
and we can apply Theorem~\ref{thmG} to $U_n^{(k)}$ and $W_n^{(\ell)}$.
By Claim~\ref{cl8}, we may write
\begin{align*}
	\mu_1(U_n)&=\sum_{k\leqslant\frac{n-t-1}{2}}\!\! |U_n^{(k)}|\phi_1(k)+ o(q^{-2tn}), \\
	\mu_2(W_n)&=\sum_{\ell\leqslant\frac{n-t-1}{2}}\!\!|W_n^{(\ell)}|\phi_2(\ell)+ o(q^{-2tn}).
\end{align*}
Then, by Theorem~\ref{thmG}, we have
\begin{align*}
\MoveEqLeft[7] \left(\sum_{k\leqslant\frac{n-t-1}{2}}\!\!|U_n^{(k)}|\phi_1(k)\right)\!\!\!
\left(\sum_{\ell\leqslant\frac{n-t-1}{2}}\!\!|W_n^{(\ell)}|\phi_2(\ell)\right) \\
&=\sum_{k\leqslant\frac{n-t-1}{2}}\sum_{\ell\leqslant\frac{n-t-1}{2}}\!\!|U_n^{(k)}||W_n^{(\ell)}|
\phi_1(k)\phi_2(\ell)\\
&\leqslant
\sum_{k\leqslant\frac{n-t-1}{2}}\sum_{\ell\leqslant\frac{n-t-1}{2}}\qbinom{n-t}{k-t}\qbinom{n-t}{\ell-t}
\phi_1(k)\phi_2(\ell)\\
&=
\left(\sum_{k\leqslant\frac{n-t-1}{2}}\qbinom{n-t}{k-t}\phi_1(k)\right) \!\!\!
\left(\sum_{\ell\leqslant\frac{n-t-1}{2}}\qbinom{n-t}{\ell-t}\phi_2(\ell)\right)
\\
&\leqslant \mu_1(A_n^{(t)})\mu_2(A_n^{(t)}),
\end{align*}
where $A_n^{(t)}$ is defined by \eqref{def:An}.
It follows that
\begin{align*}
\mu_1(U_n)\mu_2(W_n)
&= \left(\sum_{k\leqslant\frac{n-t-1}{2}}\!\!|U_n^{(k)}|\phi_1(k)\right)\!\!\!
\left(\sum_{\ell\leqslant\frac{n-t-1}{2}}\!\!|W_n^{(\ell)}|\phi_2(\ell)\right) \!+ o(q^{-2tn}) \\
&\leqslant
\mu_1(A_n^{(t)})\mu_2(A_n^{(t)})+ o(q^{-2tn}).
\end{align*}
Thus, by using Lemma~\ref{la2}, we have
\begin{equation*}
	g(n,t,\sigma_1,\sigma_2)=
	\left(\mu_1(U_n) \mu_2(W_n)\right)^{\frac{1}{n}}  
	\leqslant \left(\mu_1(A_n^{(t)})\mu_2(A_n^{(t)})+ o(q^{-2tn}) \right)^{\frac{1}{n}}
	\rightarrow q^{-(2-\theta_1-\theta_2)t}.
\end{equation*}
The opposite inequality follows from
\begin{equation*}
	g(n,t,\sigma_1,\sigma_2)\geqslant\left(\mu_1(A_n^{(t)})\mu_2(A_n^{(t)})\right)^{\frac{1}{n}}\rightarrow q^{-(2-\theta_1-\theta_2)t}.
\end{equation*}
This completes the proof of Theorem \ref{thm3}.
\end{proof}

\section{Concluding remarks}

\subsection{Theorem~\ref{thm1} from Hoffman's bound}

Here, we show that the inequality in Theorem~\ref{thm1} can also be interpreted as an application of the so-called Hoffman's bound.
Let $\mathcal{G}=(\Omega,E)$ be a finite
simple graph with vertex set $\Omega$ and edge set $E$, and let $\phi:\Omega\to[0,1]$ be a weight function such that
$\sum_{x\in\Omega}\phi(x)=1$.
Let $\mu:2^{\Omega}\to[0,1]$ be the probability measure on $\Omega$ defined by
$\mu(U):=\sum_{x\in U}\phi(x)$.
We say that $U\subset\Omega$ is independent if 
no two vertices in $U$ are adjacent in $\mathcal{G}$.
The independence measure of $\mathcal{G}$, denoted by $\alpha(\mathcal{G})$, is the
maximum of $\mu(U)$ over independent sets $U\subset\Omega$.
We say that a symmetric matrix 
$B\in{\mathbb R}^{\Omega\times\Omega}$ reflects adjacency in $\mathcal{G}$ if 
\begin{itemize}
\item $B_{x,y}=0$ whenever $\{x,y\}\not\in E$, and
\item $B$ has an eigenvector $\bm{w}\in{\mathbb R}^\Omega$ such that
$\bm{w}_x=\sqrt{\phi(x)}=\sqrt{\mu(x)}$ \,$(x\in\Omega)$.
\end{itemize}
Let $\bm{u}_1,\bm{u}_2,\dots,\bm{u}_{|\Omega|}$ be eigenvectors of
$B$. We may assume that they form an orthonormal basis with respect to the
standard inner product, and that $\bm{u}_1=\bm{w}$. Let $\lambda_i$ denote the
eigenvalue for $\bm{u}_i$, and let 
$\lambda_{\min}:=\min\{\lambda_i:2\leqslant i\leqslant |\Omega|\}$.
Under these assumptions, we have the following upper bound for $\alpha(\mathcal{G})$
in terms of $\lambda_1$ and $\lambda_{\min}$.

\begin{lemma}[Hoffman's bound]
We have $\alpha(\mathcal{G})\leqslant\frac{-\lambda_{\min}}{\lambda_1-\lambda_{\min}}$.
\end{lemma}

The original version of Hoffman's bound assumes that $\mathcal{G}$ is a regular graph and that $\mu$ is the uniform measure on $\Omega$, but its proof works for the above general version: evaluate $\bm{v}^{\mathsf{T}}B\bm{v}$ in two ways for an independent set $U\subset \Omega$, where $\bm{v}\in\mathbb{R}^{\Omega}$ is defined by
\begin{equation*}
	\bm{v}_x=\begin{cases} \sqrt{\mu(x)} & \text{if} \ x\in U, \\ 0 & \text{if} \ x\not\in U, \end{cases} \qquad (x\in \Omega).
\end{equation*}

Let the vertex set of $\mathcal{G}$ be $\Omega:=\Omega_n$, where two distinct vertices $x,y\in\Omega_n$ are adjacent whenever $x\cap y=0$.
Observe that $U\subset\Omega_n$ is independent
in $\mathcal{G}$ if and only if $U$ is an intersecting family. 
Let $\phi:=\phi_{\sigma,n}$, $\mu:=\mu_{\sigma}$,
and $B:=A'$ from \eqref{A'} with \eqref{our choice for A'}.
Then it follows that $A'$ reflects the adjacency, and that $\bm{w}=\Delta^{\frac{1}{2}}\bm{1}=\bm{v}_0$ is a $\frac{1}{1+\sigma}$-eigenvector of $A'$ (see \eqref{w0} and \eqref{v0}).
Since $\lambda_1=\frac{1}{1+\sigma}$ and
$\lambda_{\min}=-\frac{\sigma}{1+\sigma}$ by \eqref{0th block eigenvalues} and \eqref{ith block eigenvalues}, Hoffman's bound 
yields $\alpha(\mathcal{G})\leqslant\frac{\sigma}{1+\sigma}$. Moreover, in this case,
$\alpha(\mathcal{G})$ is the maximum $\mu_{\sigma}$-biased measure of intersecting
families in $\Omega_n$, and so we obtain the inequality in Theorem~\ref{thm1}.
In fact, we can apply our semidefinite programming formulation conversely to prove Hoffman's bound in general.
This was already done in \cite{Lovasz1979IEEE} for the original version.
We adopted this formulation in this paper for possible extendability; see, e.g., \cite{STT2017MP}.

\subsection{Comparison of Theorem~\ref{thmG} with the corresponding subset version}

For the case $t=1$, Theorem~\ref{thmG} reads as follows.

\begin{cor}
Let $n\geqslant k+\ell+2$. If $U\subset\ksubsp{n}{k}$ and $W\subset\ksubsp{n}{\ell}$ 
are cross $1$-intersecting, then $|U||W|\leqslant\qbinom{n-1}{k-1}\qbinom{n-1}{\ell-1}$.
\end{cor}
If we replace $\ksubsp{n}{k}$ and $\ksubsp{n}{\ell}$ in the above result with
$\binom{X_n}{k}$ and $\binom{X_n}{\ell}$, respectively, then the situation becomes
more complicated; see \cite{T2020}.
Indeed, for any large $M>0$, we can find
$n,k,\ell$ with $n=k+\ell+M$ and cross $1$-intersecting families 
$U\subset\binom{X_n}{k}$ and $W\subset\binom{X_n}{\ell}$ such that
\begin{equation}\label{UW>EKR}
	|U||W|>\binom{n-1}{k-1}\!\binom{n-1}{\ell-1}. 
\end{equation}
To see this, let
\begin{equation*}
	U:=\left\{x\in\binom{X_n}k:x\cap X_2\neq\emptyset\right\}, \qquad W:=\left\{x\in\binom{X_n}{\ell}:X_2\subset x\right\}.
\end{equation*}
Then, $U$ and $W$ are cross $1$-intersecting, and $|U|=\binom{n}{k}-\binom{n-2}{k}$, $|W|=\binom{n-2}{\ell-2}$. In this case, \eqref{UW>EKR}
is equivalent to $(2n-k-1)k(\ell-1)>(n-1)^2$,
and this condition is satisfied, for example, if $k\geqslant 2$, $\ell=9k$, 
$n=17k$. This means that for any given $M=7k$, we can construct families
satisfying \eqref{UW>EKR}, and so the condition $n\geqslant k+\ell+M$ is not
sufficient to get the upper bound $\binom{n-1}{k-1}\binom{n-1}{\ell-1}$ for
the product of the sizes of cross $1$-intersecting families.

It is interesting to determine whether the condition $n\geqslant k+\ell+2$ 
in the corollary is sharp or not. 
We cannot replace the condition with $n\geqslant k+\ell$.
To see this, let $n=k+\ell$, $k=1$, $\ell\geqslant 3$, and fix $z\in\ksubsp{n}{2}$.
Then
\begin{equation*}
	U:=\{x\in\ksubsp{n}{1}:x\subset z\}, \qquad W:=\{x\in\ksubsp{n}{\ell}:z\subset x\}
\end{equation*}
are cross $1$-intersecting, and
\begin{equation*}
	|U||W|=\qbinom{2}{1}\qbinom{n-2}{\ell-2}=\qbinom{2}{1}\qbinom{\ell-1}{1}
>\qbinom{\ell}{1}=\qbinom{n-1}{k-1}\qbinom{n-1}{\ell-1}.
\end{equation*}

\section*{Acknowledgments}

HT was supported by JSPS KAKENHI Grant Number JP23K03064.
NT was supported by JSPS KAKENHI Grant Number JP23K03201.

\end{document}